# On some new formulae involving the Stieltjes constants

Donal F. Connon

1 February 2019

**Abstract**

We show that the generalised Stieltjes constants $\gamma_n(x)$ may be represented by infinite series involving logarithmic terms. In particular we have

$$\gamma_n(x) = -\frac{\log^{n+1} x}{n+1} + \sum_{k=0}^{\infty} \left[ \frac{\log^n(k+x)}{k+x} - \frac{1}{n+1}[\log^{n+1}(k+1+x) - \log^{n+1}(k+x)] \right]$$

and for $x=1$ we have

$$\gamma_n = \sum_{k=1}^{\infty} \left[ \frac{\log^n k}{k} - \frac{1}{n+1}[\log^{n+1}(k+1) - \log^{n+1} k] \right]$$

Some relations involving the derivatives of the Hurwitz zeta function $\varsigma^{(n)}(0,x)$ are also investigated. These include

$$(-1)^{k+1}[\varsigma^{(k+1)}(0,x) - \varsigma^{(k+1)}(0)] = \log^{k+1} x + \sum_{n=1}^{\infty} \left[ \log^{k+1}(n+x) - \log^{k+1} n - x[\log^{k+1}(n+1) - \log^{k+1} n] \right]$$

**1. The Hurwitz zeta function** $\varsigma(s, x)$

The Hurwitz zeta function $\varsigma(s, x)$ is initially defined for $\operatorname{Re} s > 1$ and $x > 0$ by

(1.1) $$\varsigma(s,x) = \sum_{n=0}^{\infty} \frac{1}{(n+x)^s}$$

Note that $\varsigma(s, x)$ may be analytically continued to the whole $s$ plane except for a simple pole at $s=1$. For example, Hasse (1898-1979) showed that [26]

(1.2) $$\varsigma(s,x) = \frac{1}{s-1} \sum_{n=0}^{\infty} \frac{1}{n+1} \sum_{k=0}^{n} \binom{n}{k} \frac{(-1)^k}{(k+x)^{s-1}}$$

is a globally convergent series for $\varsigma(s, x)$ and, except for $s = 1$, provides an analytic continuation of $\varsigma(s, x)$ to the entire complex plane.

It may be noted from (1.1) that $\varsigma(s,1) = \varsigma(s)$.

We easily see from (1.2) that

$$\lim_{s \to 1}[(s-1)\varsigma(s,x)] = \sum_{n=0}^{\infty} \frac{1}{n+1} \sum_{k=0}^{n} \binom{n}{k}(-1)^k$$

and, since $(1-1)^n = \sum_{k=0}^{n} \binom{n}{k}(-1)^k = \delta_{n,0}$, we have

(1.3) $$\lim_{s \to 1}[(s-1)\varsigma(s,x)] = 1$$

which shows that $\varsigma(s,x)$ has a simple pole at $s=1$. This enables us to write the Laurent expansion for the Hurwitz zeta function shown in (2.1) below.

## 2. The generalised Stieltjes constants $\gamma_n(x)$

The generalised Stieltjes constants $\gamma_n(x)$ are the coefficients in the Laurent expansion of the Hurwitz zeta function $\varsigma(s,u)$ about $s=1$

(2.1) $$\varsigma(s,x) = \sum_{n=0}^{\infty} \frac{1}{(n+x)^s} = \frac{1}{s-1} + \sum_{n=0}^{\infty} \frac{(-1)^n}{n!} \gamma_n(x)(s-1)^n$$

We have

(2.2) $$\gamma_0(x) = -\psi(x)$$

where $\psi(x)$ is the digamma function which is the logarithmic derivative of the gamma function $\psi(x) = \frac{d}{dx}\log\Gamma(x)$. It is easily seen from the definition of the Hurwitz zeta function that $\varsigma(s,1) = \varsigma(s)$ and accordingly that $\gamma_n(1) = \gamma_n$.

The Stieltjes constants $\gamma_n$ (or Euler-Mascheroni constants) are the coefficients of the Laurent expansion of the Riemann zeta function $\varsigma(s)$ about $s=1$.

(2.3) $$\varsigma(s) = \frac{1}{s-1} + \sum_{n=0}^{\infty} \frac{(-1)^n}{n!} \gamma_n (s-1)^n$$

Since $\lim_{s \to 1}\left[\varsigma(s) - \frac{1}{s-1}\right] = \gamma$ it is clear that $\gamma_0 = \gamma$. An elementary proof of $\gamma_0(x) = -\psi(x)$ was recently given by the author in [21] (this formula was first obtained by Berndt [2] in 1972).

We note from (2.1) that

(2.4) $$\frac{\partial}{\partial s}[(s-1)\varsigma(s,x)]\bigg|_{s=1} = \gamma_0(x) = -\psi(x)$$

and it is easy to see that

$$\varsigma^{(m)}(s,x) - \varsigma^{(m)}(s,y) = \sum_{n=0}^{\infty} \frac{(-1)^n}{n!} [\gamma_n(x) - \gamma_n(y)] n(n-1)\ldots(n-m+1)(s-1)^{n-m}$$

In the limit as $s \to 1$ this becomes

(2.5) $\quad \lim_{s \to 1} [\varsigma^{(m)}(s,x) - \varsigma^{(m)}(s,y)] = (-1)^m [\gamma_m(x) - \gamma_m(y)]$

and hence using (1.1) we have for $m \geq 1$

(2.6) $\quad \gamma_m(x) - \gamma_m(y) = \sum_{n=0}^{\infty} \left[ \frac{\log^m(n+x)}{n+x} - \frac{\log^m(n+y)}{n+y} \right]$

and in particular

(2.7) $\quad \gamma_m(x) - \gamma_m = \sum_{n=0}^{\infty} \left[ \frac{\log^m(n+x)}{n+x} - \frac{\log^m(n+1)}{n+1} \right]$

and $m = 0$ gives us

(2.8) $\quad \psi(x) - \gamma = -\sum_{n=0}^{\infty} \left[ \frac{1}{n+x} - \frac{1}{n+1} \right]$

We see from (1.1) that

(2.9) $\quad \varsigma(s, 1+x) = \varsigma(s,x) - \frac{1}{x^s}$

and hence taking the limit as $s \to 1$ we have

(2.10) $\quad \gamma_m(1+x) - \gamma_m(x) = -\frac{\log^m x}{x}$

In the particular case $m = 0$ we have the familiar formula for the digamma function

$$\psi(1+x) = \psi(x) + \frac{1}{x}$$

$\square$

We have for example using (1.1) for $x > 0$ and $1 - x > 0$ (i.e. for $0 < x < 1$)

(2.11) $\quad \varsigma(s+1, x) - \varsigma(s+1, 1-x) = \sum_{n=0}^{\infty} \left[ \frac{1}{(n+x)^{s+1}} - \frac{1}{(n+1-x)^{s+1}} \right]$

and using (2.1) we have

$$\varsigma(s+1,x) - \varsigma(s+1,1-x) = \sum_{n=0}^{\infty} \frac{(-1)^n}{n!} [\gamma_n(x) - \gamma_n(1-x)] s^n$$

In the limit as $s \to 0$ we obtain

$$\gamma_0(x) - \gamma_0(1-x) = \sum_{n=0}^{\infty} \left[ \frac{1}{n+x} - \frac{1}{n+1-x} \right]$$

or equivalently

$$\psi(1-x) - \psi(x) = \frac{1}{x} + \sum_{n=1}^{\infty} \frac{1}{x^2 - n^2}$$

Since

$$\cot(\pi x) = \psi(1-x) - \psi(x)$$

we have therefore obtained a rather straightforward derivation of the decomposition formula for the cotangent function

(2.12) $$\cot(\pi x) = \frac{1}{x} + \sum_{n=1}^{\infty} \frac{1}{x^2 - n^2}$$

□

Differentiating (2.9) gives us

$$\varsigma^{(n)}(s, 1+x) = \varsigma^{(n)}(s, x) + (-1)^{n+1} \frac{\log^n x}{x^s}$$

so that with $s = 0$ we have

(2.13) $$\varsigma^{(n)}(0, 1+x) = \varsigma^{(n)}(0, x) + (-1)^{n+1} \log^n x$$

This shows that $\varsigma^{(n)}(0, 2) = \varsigma^{(n)}(0)$.

□

We see from (1.1) that

$$\frac{\partial}{\partial x} \varsigma(s, x) = -s\varsigma(s+1, x)$$

and using (2.1) we find that

(2.14) $$s\varsigma(s+1, x) = 1 + \sum_{k=0}^{\infty} \frac{(-1)^k}{k!} \gamma_k(x) s^{k+1}$$

Differentiation results in

$$\frac{\partial^{n+1}}{\partial s^{n+1}}[s\varsigma(s+1,x)] = \sum_{k=0}^{\infty}\frac{(-1)^k}{k!}\gamma_k(x)(k+1)k(k-1)...(k-n+1)s^{k-n}$$

or equivalently

$$\frac{\partial^{n+1}}{\partial s^{n+1}}\frac{\partial}{\partial x}\varsigma(s,x) = -\sum_{k=0}^{\infty}\frac{(-1)^k}{k!}\gamma_k(x)(k+1)k(k-1)...(k-n+1)s^{k-n}$$

We note that the partial derivatives commute in the region where $\varsigma(s,x)$ is analytic and hence we have

$$\frac{\partial^n}{\partial s^n}\frac{\partial}{\partial x}\varsigma(s,x) = \frac{\partial}{\partial x}\frac{\partial^n}{\partial s^n}\varsigma(s,x)$$

We therefore obtain

(2.15) $$\quad\frac{\partial}{\partial x}\varsigma^{(n+1)}(0,x) = (n+1)(-1)^{n+1}\gamma_n(x)$$

and integration results in

(2.16) $$\quad\int_1^u \gamma_n(x)\,dx = \frac{(-1)^{n+1}}{n+1}\left[\varsigma^{(n+1)}(0,u) - \varsigma^{(n+1)}(0)\right]$$

It is immediately seen that Lerch's formula [4]

(2.17) $$\quad\varsigma'(0,x) = \log\Gamma(x) + \varsigma'(0)$$

arises in the case $n=0$ because $\gamma_0(x) = -\psi(x)$.

We have for $u=2$ in (2.16)

(2.18) $$\quad\int_1^2 \gamma_n(x)\,dx = \frac{(-1)^{n+1}}{n+1}\left[\varsigma^{(n+1)}(0,2) - \varsigma^{(n+1)}(0)\right]$$

and, because $\varsigma^{(n)}(0,2) = \varsigma^{(n)}(0)$, we deduce that

(2.19) $$\quad\int_1^2 \gamma_n(x)\,dx = 0$$

As we shall see later, this may be directly confirmed by integrating (3.14).

The vanishing of the integral in (2.19) tells us that the function $\gamma_n(x)$ must assume both positive and negative values in the interval $[1,2]$ (and it is well known that this is indeed the case with the digamma function $\psi(x)$, $\gamma_0(x) = -\psi(x)$, which is negative for $x \in (0,\alpha)$ and positive for $x > \alpha$ where $\alpha \simeq 1.461632144968...$). The intermediate value theorem then dictates that $\gamma_n(x)$ must have at least one real zero in the interval $[1, 2]$.

We note from (2.10) that

(2.20) $\qquad \gamma_n(1) = \gamma_n(2)$ for $n \geq 1$

In view of this equality, we see that for $n \geq 1$ $\gamma_n(x)$ must have at least two real zeros in the interval $[1, 2]$. Furthermore, there must exist at least one $c \in [1, 2]$ such that $\gamma_n'(c) = 0$.

Using (2.10) and (2.19) we see that

$$\int_0^1 \left[ \gamma_n(x) - \frac{\log^n x}{x} \right] dx = \int_0^1 \gamma_n(1+x) \, dx$$

$$= \int_1^2 \gamma_n(t) \, dt = 0$$

Since $\log^{2n} x > 0$ for $x \in (0,1)$ we deduce that $\gamma_{2n}(x)$ must be positive in some sub-interval of $[0,1]$. Similarly, we deduce that $\gamma_{2n+1}(x)$ must be negative in some sub-interval of $[0,1]$.

We note that (2.19) may be expressed as

$$\int_0^1 \gamma_n(x-1) \, dx = 0$$

We know from (3.14) that

$$\lim_{x \to 0} \gamma_{2n}(x) = +\infty$$

$$\lim_{x \to 0} \gamma_{2n+1}(x) = -\infty$$

$\square$

We have Raabe's integral $\int_0^1 \log \Gamma(x) \, dx = \frac{1}{2} \log(2\pi)$ and we deduce from Lerch's formula (2.17) that

(2.21) $\qquad \int_0^1 \varsigma'(0, x) \, dx = 0$

$\square$

The values of $\gamma_1(u)$ in the case where $u$ is a rational number, $u = \frac{p}{q}$, are known: there are several recent relevant papers by Blagouchine [6], Coffey [17] and the author ([18] and [20]) on this topic.

(2.22) $$\gamma_1\left(\frac{p}{q}\right) = \gamma_1 + [\gamma + \log(2\pi q)]\left[\gamma + \psi\left(\frac{p}{q}\right)\right] + \sum_{j=1}^{q-1} \cos\left(\frac{2\pi jp}{q}\right)\varsigma''\left(0, \frac{j}{q}\right)$$

$$+ \pi \sum_{j=1}^{q-1} \sin\left(\frac{2\pi jp}{q}\right)\log\Gamma\left(\frac{j}{q}\right) + \frac{1}{2}\log^2 q + \log q \log(2\pi)$$

and some examples of the first Stieltjes constants are illustrated below:

$$\gamma_1\left(\frac{1}{4}\right) = \frac{1}{2}[2\gamma_1 - 7\log^2 2 - 6\gamma \log 2] - \frac{1}{2}\pi\left[\gamma + 4\log 2 + 3\log \pi - 4\log\Gamma\left(\frac{1}{4}\right)\right]$$

$$\gamma_1\left(\frac{3}{4}\right) = \frac{1}{2}[2\gamma_1 - 7\log^2 2 - 6\gamma \log 2] + \frac{1}{2}\pi\left[\gamma + 4\log 2 + 3\log \pi - 4\log\Gamma\left(\frac{1}{4}\right)\right]$$

$$\gamma_1\left(\frac{1}{5}\right) = \frac{1}{4}\left[4\gamma_1 - \frac{5}{2}\log^2 5 - 5\gamma \log 5\right] - \frac{1}{2}\pi[\log(10\pi) + \gamma]\cot\left(\frac{\pi}{5}\right)$$

$$+ \frac{1}{4}\sqrt{5}\left[\varsigma''\left(0,\frac{1}{5}\right) - \varsigma''\left(0,\frac{2}{5}\right) - \varsigma''\left(0,\frac{3}{5}\right) + \varsigma''\left(0,\frac{4}{5}\right) - [\gamma + \log(2\pi)]\log\frac{1}{2}(3+\sqrt{5})\right]$$

$$+ \pi\left\{\sin\left(\frac{2\pi}{5}\right)\left[\log\Gamma\left(\frac{1}{5}\right) - \log\Gamma\left(\frac{4}{5}\right)\right] + \sin\left(\frac{\pi}{5}\right)\left[\log\Gamma\left(\frac{2}{5}\right) - \log\Gamma\left(\frac{3}{5}\right)\right]\right\}$$

We also have similar expressions for $\gamma_0\left(\frac{p}{q}\right)$, namely

$$\psi\left(\frac{p}{q}\right) = -\gamma - \log(2\pi q) - \frac{\pi}{2}\cot\left(\frac{p\pi}{q}\right) - 2\sum_{r=1}^{q}\cos\left(\frac{2\pi rp}{q}\right)\log\Gamma\left(\frac{r}{q}\right)$$

The values of $\gamma_n(u)$ are also known [20] in the case where $u$ is a rational number.

### 3. An infinite series for the generalised Stieltjes constants $\gamma_n(x)$

It is well known that the Stieltjes constants $\gamma_n$ may be represented by the limit

(3.1) $$\gamma_n = \lim_{N \to \infty} \sum_{k=1}^{N}\left[\frac{\log^n k}{k} - \frac{\log^{n+1} N}{n+1}\right]$$

or equivalently by

(3.2) $$\gamma_n = \lim_{N \to \infty} \sum_{k=1}^{N}\left[\frac{\log^n k}{k} - \int_{1}^{N}\frac{\log^n x}{x}dx\right]$$

Blagouchine's paper [6] contains a detailed list of references for the generalised Stieltjes constants $\gamma_n(x)$.

Osler [34)] has provided an interesting derivation of the above formula using partial sums of the Riemann zeta function.

In 1972 Berndt [2] showed that

$$(3.3) \qquad \gamma_n(x) = \lim_{N \to \infty} \left[ \sum_{k=0}^{N} \frac{\log^n(k+x)}{k+x} - \frac{\log^{n+1}(N+x)}{n+1} \right]$$

A new proof of (3.3) was recently given by the author [19].

Simple algebra gives us

$$\sum_{k=0}^{N} \frac{\log^n(k+1)}{k+1} - \frac{\log^{n+1}(N+1)}{n+1} = \sum_{k=1}^{N} \frac{\log^n k}{k} + \frac{\log^n(N+1)}{N+1} - \frac{\log^{n+1}(N+1)}{n+1}$$

$$= \sum_{k=1}^{N} \frac{\log^n k}{k} - \frac{\log^{n+1} N}{k+1} + \frac{1}{n+1} \left[ \log^{n+1} N - \log^{n+1}(N+1) \right] + \frac{\log^n(N+1)}{N+1}$$

Letting $f(x) = \log^{n+1}(N+x)$, the mean value theorem for derivatives tells us that

$$\log^{n+1}(N+x) - \log^{n+1} N = \frac{(n+1)x \log^n(N+\theta_n)}{N+\theta_n} \quad \text{where } 0 < \theta_n < x$$

and it is easily shown that this expression vanishes as $N \to \infty$. Indeed, successive applications of L'Hôpital's rule readily shows that both $\lim\limits_{N \to \infty} \frac{\log^n(N+\theta_n)}{N+\theta_n} = 0$ and $\lim\limits_{N \to \infty} \frac{\log^n(N+1)}{N+1} = 0$ and hence we have

$$(3.4) \qquad \lim_{N \to \infty} \left[ \sum_{n=1}^{N} \frac{\log^n(k+1)}{k+1} - \frac{\log^{n+1}(N+1)}{n+1} \right] = \lim_{N \to \infty} \sum_{k=1}^{N} \left[ \frac{\log^n k}{k} - \frac{\log^{n+1} N}{n+1} \right]$$

Using (3.4) we see that Stieltjes's formula (3.1) for $\gamma_n$ immediately follows by letting $x=1$ in (3.3).

Dilcher [24] also showed by a different method that $\lim\limits_{N \to \infty}[\log^{n+1}(N+x) - \log^{n+1} N] = 0$.

We shall find that the following lemma is exceedingly useful.

**Lemma 3.1**

$$(3.5) \quad \log^{n+1}(N+x) = \log^{n+1}(1+x) - \int_N^{N+1} \frac{\log^n(x+t)}{x+t} dt + \sum_{k=1}^N [\log^{n+1}(k+1+x) - \log^{n+1}(k+x)]$$

With $x = 0$ this becomes

$$(3.6) \quad \log^{n+1} N = -\int_N^{N+1} \frac{\log^n t}{t} dt + \sum_{k=1}^N [\log^{n+1}(k+1) - \log^{n+1} k]$$

where $N \geq 1$ is a positive integer.

In particular, note that the single logarithmic term $\log^{n+1}(N+x)$ has been replaced by a finite sum together with a residual term that vanishes as $N \to \infty$. This form assumes importance in various applications encountered in this paper.

**Proof**

It is easily seen that

$$\int_1^N \frac{\log^n(x+t)}{x+t} dt = \frac{1}{n+1} \left[ \log^{n+1}(N+x) - \log^{n+1}(1+x) \right]$$

so that

$$\log^{n+1}(N+x) = \log^{n+1}(1+x) + (n+1) \int_1^N \frac{\log^n(x+t)}{x+t} dt$$

We then create the requisite summation by splitting the integral as follows

$$\int_1^N \frac{\log^n(x+t)}{x+t} dt = \sum_{k=1}^{N-1} \int_k^{k+1} \frac{\log^n(x+t)}{x+t} dt$$

$$= \sum_{k=1}^N \int_k^{k+1} \frac{\log^n(x+t)}{x+t} dt - \int_N^{N+1} \frac{\log^n(x+t)}{x+t} dt$$

$$= \frac{1}{n+1} \sum_{k=1}^N [\log^{n+1}(k+1+x) - \log^{n+1}(k+x)] - \int_N^{N+1} \frac{\log^n(x+t)}{x+t} dt$$

and hence we have

$$\log^{n+1}(N+x) = \log^{n+1}(1+x) - \int_N^{N+1} \frac{\log^n(x+t)}{x+t} dt + \sum_{k=1}^N [\log^{n+1}(k+1+x) - \log^{n+1}(k+x)]$$

This proves Lemma 3.1. Whilst this derivation superficially involve integration, it should be noted that it is purely algebraic in character (i.e. it is a finite telescoping sum).

The above identities may be represented in a slightly different manner as follows.

We have
$$\varsigma(s, 1+x) - \varsigma(s, x) = -\frac{1}{x^s}$$

and differentiation with respect to $s$ gives us

$$\varsigma^{(n)}(s, 1+x) - \varsigma^{(n)}(s, x) = \frac{(-1)^{n+1} \log^n x}{x^s}$$

Letting $s = 0$ results in

$$\varsigma^{(n+1)}(0, 1+x) - \varsigma^{(n+1)}(0, x) = (-1)^n \log^{n+1} x$$

and this implies that

$$\varsigma^{(n+1)}(0, 1+k+1+x) - \varsigma^{(n+1)}(0, k+1+x) = (-1)^n \log^{n+1}(k+1+x)$$

$$\varsigma^{(n+1)}(0, 1+k+x) - \varsigma^{(n+1)}(0, k+x) = (-1)^n \log^{n+1}(k+x)$$

We then obtain

$$\varsigma^{(n+1)}(0, k+2+x) - 2\varsigma^{(n+1)}(0, k+1+x) + \varsigma^{(n+1)}(0, k+x)$$

$$= (-1)^n [\log^{n+1}(k+1+x) - \log^{n+1}(k+x)]$$

Therefore (3.5) may be written as

**Lemma 3.2**

(3.7) $$\log^{n+1}(N+x) = (-1)^n \sum_{k=1}^{N} \left\{ \varsigma^{(n+1)}(0, k+2+x) - 2\varsigma^{(n+1)}(0, k+1+x) + \varsigma^{(n+1)}(0, k+x) \right\}$$

$$+ \log^{n+1}(1+x) - \int_{N}^{N+1} \frac{\log^n(x+t)}{x+t} dt$$

With $x = 0$ this becomes

(3.8) $$\log^{n+1} N = (-1)^n \sum_{k=1}^{N} \left\{ \varsigma^{(n+1)}(0, k+2) - 2\varsigma^{(n+1)}(0, k+1) + \varsigma^{(n+1)}(0, k) \right\} - \int_{N}^{N+1} \frac{\log^n t}{t} dt$$

where $N \geq 1$ is a positive integer. These expressions at least highlight the intimate connection with the Hurwitz zeta function.

**Proposition 3.1**

We have for $x > 0$

(3.9) $$\gamma_n(x) = -\frac{\log^{n+1} x}{n+1} + \sum_{k=0}^{\infty}\left[\frac{\log^n(k+x)}{k+x} - \frac{1}{n+1}[\log^{n+1}(k+1+x) - \log^{n+1}(k+x)]\right]$$

and

(3.10) $$\gamma_n = \sum_{k=1}^{\infty}\left[\frac{\log^n k}{k} - \frac{1}{n+1}[\log^{n+1}(k+1) - \log^{n+1} k]\right]$$

**Proof**

As mentioned previously in (3.3), it is well known that for $x > 0$

$$\gamma_n(x) = \lim_{N\to\infty}\left[\sum_{k=0}^{N}\frac{\log^n(k+x)}{k+x} - \frac{\log^{n+1}(N+x)}{n+1}\right]$$

and using Lemma 3.1

$$\log^{n+1}(N+x) = \log^{n+1}(1+x) - \int_N^{N+1}\frac{\log^n(x+t)}{x+t}dt + \sum_{k=1}^{N}[\log^{n+1}(k+1+x) - \log^{n+1}(k+x)]$$

we deduce that

$$\sum_{k=0}^{N}\frac{\log^n(k+x)}{k+x} - \frac{\log^{n+1}(N+x)}{n+1}$$

$$= \frac{\log^n x}{x} + \sum_{k=1}^{N}\left[\frac{\log^n(k+x)}{k+x} - \frac{1}{n+1}\left\{\log^{n+1}(k+1+x) - \log^{n+1}(k+x)\right\}\right]$$

$$-\frac{1}{n+1}\log^{n+1}(1+x) + \frac{1}{n+1}\int_N^{N+1}\frac{\log^n(x+t)}{x+t}dt$$

where we have now started the summation at $k=1$.

Therefore, as $N \to \infty$ we obtain

(3.11) $$\gamma_n(x) = \frac{\log^n x}{x} - \frac{\log^{n+1}(1+x)}{n+1} + \sum_{k=1}^{\infty}\left\{\frac{\log^n(k+x)}{k+x} - \frac{1}{n+1}[\log^{n+1}(k+1+x) - \log^{n+1}(k+x)]\right\}$$

since $\lim_{N\to\infty}\int_N^{N+1}\frac{\log^n(x+t)}{x+t}dt = 0$. This limit appears to be obvious but a separate evaluation is outlined below.

We have

$$\int_N^{N+1} \frac{\log^n(x+t)}{x+t}\,dt = \frac{1}{n+1}[\log^{n+1}(N+x+1) - \log^{n+1}(N+x)]$$

and, as noted previously, it is easily shown that $\lim_{N\to\infty}[\log^{n+1}(N+x+1) - \log^{n+1}(N+x)] = 0$.

We may write (3.11) as

$$(3.12) \quad \gamma_n(x) = -\frac{\log^{n+1} x}{n+1} + \sum_{k=0}^{\infty}\left\{\frac{\log^n(k+x)}{k+x} - \frac{1}{n+1}[\log^{n+1}(k+1+x) - \log^{n+1}(k+x)]\right\}$$

and we see that

$$(3.13) \quad \gamma_n = \sum_{k=1}^{\infty}\left[\frac{\log^n k}{k} - \frac{\log^{n+1}(k+1) - \log^{n+1} k}{n+1}\right]$$

We shall see another representation later in (5.21)

$$(3.14) \quad \gamma_n(x) = \sum_{k=0}^{\infty}\left[\frac{\log^n(k+x)}{k+x} - \frac{1}{n+1}[\log^{n+1}(k+2) - \log^{n+1}(k+1)]\right]$$

where the variable $x$ only appears in the first term. Letting $x = 1$ in (3.14) immediately results in (3.13).

Other forms may be obtained by using (2.20): $\gamma_n(1) = \gamma_n(2)$ for $n \geq 1$.

Subtracting (3.14) from (3.12) gives us

$$\log^{n+1} x = \sum_{k=0}^{\infty}\{[\log^{n+1}(k+2) - \log^{n+1}(k+1)] - [\log^{n+1}(k+1+x) - \log^{n+1}(k+x)]\}$$

The right-hand side may be written as

$$\lim_{N\to\infty}\sum_{k=0}^{N}\{[\log^{n+1}(k+2) - \log^{n+1}(k+1)] - [\log^{n+1}(k+1+x) - \log^{n+1}(k+x)]\}$$

and using the finite telescoping sum $\sum_{k=0}^{N}[a_{k+1} - a_k] = a_{N+1} - a_0$, this becomes

$$= \lim_{N\to\infty}[\log^{n+1}(N+2) - \log^{n+1}(N+1+x) + \log^{n+1} x]$$

This confirms that (3.14) and (3.12) are equivalent.

We also have

(3.15) $$\gamma_n(x) = \frac{\log^n x}{x} - \frac{\log^{n+1}(1+x)}{n+1}$$

$$+ \sum_{k=1}^{\infty} \left\{ \frac{\log^n(k+x)}{k+x} - \frac{(-1)^n}{n+1} [\varsigma^{(n+1)}(0, k+2+x) - 2\varsigma^{(n+1)}(0, k+1+x) + \varsigma^{(n+1)}(0, k+x)] \right\}$$

(3.15.1) $$\gamma_n(x) = \frac{\log^n x}{x} - \frac{\log^{n+1}(1+x)}{n+1}$$

$$+ \sum_{k=1}^{\infty} \left\{ \frac{\log^n(k+x)}{k+x} - \frac{1}{n+1} [\log^{n+1}(k+1+x) - \log^{n+1}(k+x)] \right\}$$

Coffey [17] provided a more complicated series for $\gamma_n(x)$ in 2011

(3.16) $$\gamma_n(x) = \sum_{k=0}^{m} \frac{\log^n(k+x)}{k+x} - \frac{\log^{n+1}(m+x)}{n+1} - \frac{\log^n(m+x)}{2(m+x)}$$

$$\sum_{j=m}^{\infty} \left[ \begin{array}{c} \log^n(j+1+x) - \log^n(j+x) - \frac{1}{n+1}[\log^{n+1}(j+1+x) - \log^{n+1}(j+x)] \\ -(x+j+\tfrac{1}{2})\left(\{\Gamma[n,\log(j+x)] - \Gamma[n,\log(j+1+x)]\} - \{\Gamma[n+1,\log(j+x)] - \Gamma[n+1,\log(j+1+x)]\}\right) \end{array} \right]$$

where $\Gamma(s,t)$ is the incomplete gamma function.

With $m=0$ we have

(3.17) $$\gamma_n(x) = \frac{1}{2} \frac{\log^n x}{x} - \frac{\log^{n+1} x}{n+1}$$

$$\sum_{j=0}^{\infty} \left[ \begin{array}{c} \log^n(j+1+x) - \log^n(j+x) - \frac{1}{n+1}[\log^{n+1}(j+1+x) - \log^{n+1}(j+x)] \\ -(x+j+\tfrac{1}{2})\left(\{\Gamma[n,\log(j+x)] - \Gamma[n,\log(j+1+x)]\} - \{\Gamma[n+1,\log(j+x)] - \Gamma[n+1,\log(j+1+x)]\}\right) \end{array} \right]$$

□

As a check on our algebra, let us consider the simplest case, namely $\gamma_0$. We have using (3.15)

$$\gamma_0 = \sum_{k=1}^{\infty} \left[ \frac{1}{k} - \{\varsigma'(0, k+2) - 2\varsigma'(0, k+1) + \varsigma'(0, k)\} \right]$$

and using Lerch's formula (2.17)

$$\varsigma'(0,x) = \log\Gamma(x) - \frac{1}{2}\log(2\pi)$$

we obtain

$$\gamma_0 = \sum_{k=1}^{\infty}\left[\frac{1}{k} - \{\log\Gamma(k+2) - 2\log\Gamma(k+1) + \log\Gamma(k)\}\right]$$

It is easily seen that

$$\log\Gamma(k+2) - 2\log\Gamma(k+1) + \log\Gamma(k) = \log\frac{\Gamma(k+2)\Gamma(k)}{\Gamma(k+1)^2}$$

$$= \log\frac{(k+1)!(k-1)!}{[(k)!]^2}$$

$$= \log\left(1+\frac{1}{k}\right)$$

and we therefore obtain

$$\gamma_0 = \gamma = \sum_{k=1}^{\infty}\left[\frac{1}{k} - \log\left(1+\frac{1}{k}\right)\right]$$

□

We now briefly digress to have a leisurely review of some very familiar elementary analysis.

We have the Weierstrass canonical form of the gamma function [36, p.1]

(3.18)  $$\frac{1}{\Gamma(x)} = xe^{\gamma x}\prod_{n=1}^{\infty}\left\{\left(1+\frac{x}{n}\right)e^{-\frac{x}{n}}\right\}$$

where $\gamma$ is Euler's constant defined by $\gamma = \lim_{n\to\infty}[H_n - \log n]$ and $H_n$ are the harmonic numbers $H_n = \sum_{k=1}^{n}\frac{1}{k}$.

Taking logarithms results in

(3.19)  $$\log\Gamma(1+x) = -\gamma x - \sum_{n=1}^{\infty}\left[\log\left(1+\frac{x}{n}\right) - \frac{x}{n}\right]$$

Euler's constant $\gamma$ may be defined as the limit of the sequence

$$\gamma = \lim_{n\to\infty}[H_n - \log n]$$

$$= \lim_{n\to\infty}[H_n - \log(n+1) + \log(n+1) - \log n]$$

$$= \lim_{n\to\infty}[H_n - \log(n+1)] + \lim_{n\to\infty}\log\left(1+\frac{1}{n}\right)$$

$$= \lim_{n\to\infty}[H_n - \log(n+1)]$$

We have
$$\sum_{k=1}^{n}\log\left(1+\frac{1}{k}\right) = \sum_{k=1}^{n}[\log(k+1) - \log k]$$

and it is easily seen that this telescopes to

$$\sum_{k=1}^{n}\log\left(1+\frac{1}{k}\right) = \log(n+1)$$

We can therefore write Euler's constant $\gamma$ as the infinite series

(3.20) $$\gamma = \sum_{n=1}^{\infty}\left[\frac{1}{n} - \log\left(1+\frac{1}{n}\right)\right]$$

which we have previously encountered above.

The combination of (3.19) with (3.20) results in

$$\log\Gamma(1+x) = \sum_{n=1}^{\infty}\left[x\log\left(1+\frac{1}{n}\right) - \log\left(1+\frac{x}{n}\right)\right]$$

and differentiation gives us (an infrequently used formula)

(3.21) $$\psi(1+x) = \sum_{n=1}^{\infty}\left[\log\left(1+\frac{1}{n}\right) - \frac{1}{n+x}\right]$$

which is equivalent to the well-known representation [36, p.14]

$$\psi(x) = \lim_{N\to\infty}\left[\log N - \sum_{n=0}^{N}\frac{1}{n+x}\right]$$

It is also reported in [31, p.81] that

(3.22) $$\psi(x) = \lim_{N\to\infty}\left[\log(N+x) - \sum_{n=0}^{N}\frac{1}{n+x}\right]$$

and it is obvious that these two representations are equivalent because $\lim_{N\to\infty}[\log(N+x) - \log N] = 0$. An elementary application of Lemma 3.1 easily converts (3.22) into (3.21). □

With $x = 1$ in (3.15) we have

$$\gamma_n = -\frac{\log^{n+1} 2}{n+1} + \sum_{k=1}^{\infty} \left\{ \frac{\log^n(k+1)}{k+1} - \frac{1}{n+1}[\log^{n+1}(k+2) - \log^{n+1}(k+1)] \right\}$$

which may be written as

(3.23) $$\gamma_n = \sum_{k=0}^{\infty} \left\{ \frac{\log^n(k+1)}{k+1} - \frac{1}{n+1}[\log^{n+1}(k+2) - \log^{n+1}(k+1)] \right\}$$

Letting $n = 0$ in (3.15) gives us

$$\gamma_0(x) = -\psi(x) = \frac{1}{x} - \log(1+x) + \sum_{k=1}^{\infty} \left\{ \frac{1}{k+x} - [\log(k+1+x) - \log(k+x)] \right\}$$

so that

$$\psi(x) = -\frac{1}{x} + \log(1+x) - \sum_{k=1}^{\infty} \left\{ \frac{1}{k+x} - \log\left(1 + \frac{1}{k+x}\right) \right\}$$

or equivalently

(3.24) $$\psi(1+x) = \log(1+x) - \sum_{k=1}^{\infty} \left\{ \frac{1}{k+x} - \log\left(1 + \frac{1}{k+x}\right) \right\}$$

This formula for the digamma function concurs with [36, p.14].

With $x = 0$ in (3.24) we obtain with $\psi(1) = -\gamma$

$$\gamma = \sum_{k=1}^{\infty} \left[ \frac{1}{k} - \log\left(1 + \frac{1}{k}\right) \right]$$

We have for example

$$\gamma_1 = \sum_{k=0}^{\infty} \left\{ \frac{\log(k+1)}{k+1} - \frac{1}{2}[\log^2(k+2) - \log^2(k+1)] \right\}$$

and this logarithmic series contains some complexity because it is known that

(3.25) $$\gamma_1 = \varsigma''(0) + [\gamma + \log(2\pi)]\log(2\pi) - \frac{1}{2}\left( [\gamma + \log(2\pi)]^2 - \frac{1}{2}\varsigma(2) \right)$$

which corresponds with the result previously obtained by Ramanujan [4] and Apostol [1]

$$\varsigma''(0) = \gamma_1 + \frac{1}{2}\gamma^2 - \frac{1}{24}\pi^2 - \frac{1}{2}\log^2(2\pi)$$

Referring to (3.14) we see that

$$\gamma_1(x) = \sum_{k=0}^{\infty} \left\{ \frac{\log(k+x)}{k+x} - \frac{1}{2}[\log^2(k+2) - \log^2(k+1)] \right\}$$

In theory, one could use this to deduce the Fourier series for $\gamma_1(x)$ but the requisite integrals appear to be rather intractable.

□

Referring to (3.11) we have

$$\gamma_n(1+x) - \gamma_n(x) = \frac{\log^n(1+x)}{1+x} - \frac{\log^n x}{x} - \frac{1}{n+1}[\log^{n+1}(2+x) - \log^{n+1}(1+x)]$$

$$+ \sum_{k=1}^{\infty} \left\{ \frac{\log^n(k+1+x)}{k+1+x} - \frac{\log^n(k+x)}{k+x} - \frac{1}{n+1}[\log^{n+1}(k+2+x) - \log^{n+1}(k+1+x) - \log^{n+1}(k+1+x) + \log^{n+1}(k+x)] \right\}$$

Using the finite telescoping sum

$$\sum_{k=1}^{N}[a_{k+1} - a_k] = a_{N+1} - a_1$$

we see that the finite sum

$$\sum_{k=1}^{N} \left\{ \begin{array}{l} \frac{\log^n(k+1+x)}{k+1+x} - \frac{\log^n(k+x)}{k+x} \\ -\frac{1}{n+1}[\log^{n+1}(k+2+x) - \log^{n+1}(k+1+x) - \log^{n+1}(k+1+x) + \log^{n+1}(k+x)] \end{array} \right\}$$

$$= \frac{\log^n(N+1+x)}{N+1+x} - \frac{\log^n(1+x)}{1+x}$$

$$- \frac{1}{n+1}[\log^{n+1}(N+2+x) - \log^{n+1}(2+x) - \log^{n+1}(N+1+x) + \log^{n+1}(1+x)]$$

We then have the limit

$$\lim_{N \to \infty} \sum_{k=1}^{N} \left\{ \begin{array}{l} \frac{\log^n(k+1+x)}{k+1+x} - \frac{\log^n(k+x)}{k+x} \\ -\frac{1}{n+1}[\log^{n+1}(k+2+x) - \log^{n+1}(k+1+x) - \log^{n+1}(k+1+x) + \log^{n+1}(k+x)] \end{array} \right\}$$

$$= -\frac{\log^n(1+x)}{1+x} - \frac{1}{n+1}[-\log^{n+1}(2+x) + \log^{n+1}(1+x)]$$

because $\lim_{N\to\infty} \dfrac{\log^n(N+1+x)}{N+1+x} = 0$ and $\lim_{N\to\infty}[\log^{n+1}(N+2+x) - \log^{n+1}(N+1+x)] = 0$.

Hence, we obtain the familiar result

$$\gamma_n(1+x) - \gamma_n(x) = -\frac{\log^n x}{x}$$

which serves as a useful check on the analysis.

$\square$

Referring this time to the limit in (3.3)

$$\gamma_n(x) = \lim_{N\to\infty}\left[\sum_{k=0}^{N}\frac{\log^n(k+x)}{k+x} - \frac{\log^{n+1}(N+x)}{n+1}\right]$$

we see that

$$\gamma_n(1+x) - \gamma_n(x) = \sum_{k=0}^{\infty}\left[\frac{\log^n(k+1+x)}{k+1+x} - \frac{\log^n(k+x)}{k+x}\right]$$

Using the finite telescoping sum $\sum_{k=0}^{N}[a_{k+1} - a_k] = a_{N+1} - a_0$ and the fact that

$\lim_{N\to\infty}\dfrac{\log^n(N+1+x)}{N+1+x} = 0$ we again obtain $\gamma_n(1+x) - \gamma_n(x) = -\dfrac{\log^n x}{x}$.

$\square$

In a similar manner, we may compute the difference

$$\gamma_n(x) - \gamma_n(y) = \frac{\log^n x}{x} - \frac{\log^n y}{y} + \sum_{k=1}^{\infty}\left\{\frac{\log^n(k+x)}{k+x} - \frac{\log^n(k+y)}{k+y}\right\}$$

$$= \sum_{k=0}^{\infty}\left\{\frac{\log^n(k+x)}{k+x} - \frac{\log^n(k+y)}{k+y}\right\}$$

$$= (-1)^n \lim_{s\to 1}[\varsigma^{(n)}(s,x) - \varsigma^{(n)}(s,y)]$$

$$= \gamma_n(x) - \gamma_n(y)$$

$\square$

We write (3.11) as

$$\gamma_n(x) = \frac{\log^n x}{x} - \frac{\log^{n+1}(1+x)}{n+1} + \gamma_n^*(x)$$

where

$$\gamma_n^*(x) := \sum_{k=1}^{\infty} \left\{ \frac{\log^n(k+x)}{k+x} - \frac{1}{n+1}[\log^{n+1}(k+1+x) - \log^{n+1}(k+x)] \right\}$$

Thus (2.1) may be expressed as

$$\varsigma(s,x) = \frac{1}{s-1} + \sum_{n=0}^{\infty} \frac{(-1)^n}{n!} \left[ \frac{\log^n x}{x} - \frac{\log^{n+1}(1+x)}{n+1} \right] (s-1)^n + \sum_{n=0}^{\infty} \frac{(-1)^n}{n!} \gamma_n^*(x)(s-1)^n$$

We see that

$$\sum_{n=0}^{\infty} \frac{(-1)^n}{n!} \frac{\log^n x}{x} (s-1)^n = \frac{1}{x} \exp[-(s-1)\log x] = \frac{1}{x^s}$$

and

$$-\sum_{n=0}^{\infty} \frac{(-1)^n}{n!} \frac{\log^{n+1}(1+x)}{n+1} (s-1)^n = \frac{1}{s-1} \sum_{n=0}^{\infty} \frac{(-1)^{n+1}(s-1)^{n+1} \log^{n+1}(1+x)}{(n+1)!}$$

$$= \frac{1}{s-1}[\exp[-(s-1)\log(1+x)] - 1]$$

$$= \frac{1}{s-1}[(1+x)^{1-s} - 1]$$

Hence we obtain

$$\varsigma(s,x) = \frac{1}{x^s} + \frac{(1+x)^{1-s}}{s-1} + \sum_{n=0}^{\infty} \frac{(-1)^n}{n!} \gamma_n^*(x)(s-1)^n$$

or equivalently

$$\sum_{n=1}^{\infty} \frac{1}{(n+x)^s} = \frac{(1+x)^{1-s}}{s-1} + \sum_{n=0}^{\infty} \frac{(-1)^n}{n!} \gamma_n^*(x)(s-1)^n$$

Therefore we have

$$\gamma_n = \gamma_n^*(0) = \sum_{k=1}^{\infty} \left\{ \frac{\log^n k}{k} - \frac{1}{n+1}[\log^{n+1}(k+1) - \log^{n+1} k] \right\}$$

□

Rather belatedly I recalled that the series for $\gamma_n$ in (3.10) was not new because I subsequently remembered that it had been proved by Bohman and Fröberg [7] in 1988. Their proof is outlined below.

They noted that

$$(s-1)\varsigma(s) = \sum_{k=1}^{\infty} \frac{(s-1)}{k^s} \quad \text{and} \quad \sum_{k=1}^{\infty}\left[\frac{1}{k^{s-1}} - \frac{1}{(k+1)^{s-1}}\right] = 1$$

and, assuming that $s$ is real and greater than 1, the above two equations may be subtracted to give

(3.26) $$(s-1)\varsigma(s) = 1 + \sum_{k=1}^{\infty}\left[\frac{1}{(k+1)^{s-1}} - \frac{1}{k^{s-1}} + \frac{(s-1)}{k^s}\right]$$

The well-known limit (1.3) $\lim_{s \to 1}(s-1)\varsigma(s) = 1$ may be immediately derived from the above.

Equation (3.26) may be written as

$$(s-1)\varsigma(s) = 1 + \sum_{k=1}^{\infty}\left[\exp(-(s-1)\log(k+1)) - \exp(-(s-1)\log k) + (s-1)k^{-1}\exp(-(s-1)\log k)\right]$$

$$= 1 + \sum_{k=1}^{\infty}\left[\sum_{n=0}^{\infty}\frac{(-1)^n(s-1)^n}{n!}\left(\log^n(k+1) - \log^n k\right) + \frac{s-1}{k}\sum_{n=0}^{\infty}\frac{(-1)^n(s-1)^n}{n!}\log^n k\right]$$

Dividing by $(s-1)$ we get

$$\varsigma(s) = \frac{1}{s-1} + \sum_{n=0}^{\infty}\frac{(-1)^n}{n!}\gamma_n(s-1)^n$$

where

$$\gamma_n = \sum_{k=1}^{\infty}\left[\frac{\log^n k}{k} - \frac{\log^{n+1}(k+1) - \log^{n+1} k}{n+1}\right] = \sum_{k=1}^{\infty}\left[\frac{\log^n k}{k} - \int_k^{k+1}\frac{\log^n t}{t}dt\right]$$

which is equivalent to (3.10).

In fact, serendipitously, through a private communication received from Jacques Gélinas in January 2019, I found out that (3.10) had already been discovered by Jensen [29] in 1897.

We now bring this Section to a close by presenting a very short and direct (albeit heuristic) proof of (3.9).

We consider the integral

$$\int_0^1 \sum_{k=0}^{\infty}\left[\frac{\log^n(k+x)}{k+x} - \frac{\log^n(k+u+x)}{k+u+x}\right]du$$

and using (2.6) this becomes

$$= \int_0^1 [\gamma_n(x) - \gamma_n(u+x)]du$$

$$= \gamma_n(x) - \int_0^1 \gamma_n(u+x)\,du$$

We note that

$$I = \int_0^1 \gamma_n(u+x)\,du = \int_x^{1+x} \gamma_n(z)\,dz$$

$$= \int_1^{1+x} \gamma_n(z)\,dz - \int_1^x \gamma_n(z)\,dz$$

Using (2.16)

$$\int_1^u \gamma_n(x)\,dx = \frac{(-1)^{n+1}}{n+1}\left[\varsigma^{(n+1)}(0,u) - \varsigma^{(n+1)}(0)\right]$$

we obtain

$$I = \frac{(-1)^{n+1}}{n+1}\left[\varsigma^{(n+1)}(0,1+x) - \varsigma^{(n+1)}(0,x)\right]$$

Referring to (2.13)

$$\varsigma^{(n)}(0,1+x) = \varsigma^{(n)}(0,x) + (-1)^{n+1}\log^n x$$

then confirms that

$$I = -\frac{\log^{n+1} x}{n+1}$$

Therefore, we have

$$\int_0^1 \sum_{k=0}^{\infty}\left[\frac{\log^n(k+x)}{k+x} - \frac{\log^n(k+u+x)}{k+u+x}\right]du = \gamma_n(x) + \frac{\log^{n+1}(1+x)}{n+1}$$

We now assume that we may change the order of summation and integration

$$\sum_{k=0}^{\infty}\int_0^1\left[\frac{\log^n(k+x)}{k+x} - \frac{\log^n(k+u+x)}{k+u+x}\right]du = \int_0^1\sum_{k=0}^{\infty}\left[\frac{\log^n(k+x)}{k+x} - \frac{\log^n(k+u+x)}{k+u+x}\right]du$$

and then obtain

(3.27) $$\gamma_n(x) = -\frac{\log^{n+1}(1+x)}{n+1} + \sum_{k=0}^{\infty}\int_0^1\left[\frac{\log^n(k+x)}{k+x} - \frac{\log^n(k+u+x)}{k+u+x}\right]du$$

which corresponds with (3.9)

$$\gamma_n(x) = -\frac{\log^{n+1} x}{n+1} + \sum_{k=0}^{\infty}\left[\frac{\log^n(k+x)}{k+x} - \frac{1}{n+1}[\log^{n+1}(k+1+x) - \log^{n+1}(k+x)]\right]$$

Equation (3.27) merits further attention because, as mentioned in [17] and [18], the integral $\int_0^1 \left[\frac{\log^n(k+x)}{k+x} - \frac{\log^n(k+u+x)}{k+u+x}\right] du$ may be subject to further analysis involving the harmonic numbers and the Stirling numbers of the first kind.

## 4. Some other applications of Lemma 3.1

### (i) The gamma function

We have the Gaussian form of the gamma function [36, p.2]

$$\Gamma(x) = \lim_{N \to \infty} N^x \frac{N!}{x(1+x)\ldots(N+x)}$$

or equivalently

$$\log \Gamma(x) = \lim_{N \to \infty}\left[x \log N + \log N! - \sum_{k=0}^{N} \log(k+x)\right]$$

(4.1) $$= \lim_{N \to \infty}\left[x \log N - \log x + \sum_{k=1}^{N} \log k - \sum_{k=1}^{N} \log(k+x)\right]$$

and since $\Gamma(1+x) = x\Gamma(x)$ we have

$$\log \Gamma(1+x) = \lim_{N \to \infty}\left[x \log N + \sum_{k=1}^{N} \log k - \sum_{k=1}^{N} \log(k+x)\right]$$

We now employ Lemma 3.1

$$\log N = -\int_N^{N+1} \frac{1}{t} dt + \sum_{k=1}^{N}[\log(k+1) - \log k]$$

to obtain

$$x \log N + \sum_{k=1}^{N} \log k - \sum_{k=1}^{N} \log(k+x) = -x \int_N^{N+1} \frac{1}{t} dt + x \sum_{k=1}^{N}[\log(k+1) - \log k] - \sum_{k=1}^{N}[\log(k+x) - \log k]$$

Therefore, in a rather direct manner we obtain the well-known result [36, p.2]

(4.2) $$\log \Gamma(x+1) = \sum_{k=1}^{\infty}\left[x \log\left(1 + \frac{1}{k}\right) - \log\left(1 + \frac{x}{k}\right)\right]$$

which we have seen before in Section 3.

**(ii)** Kanemitsu et al. [40] showed in 2004 that (see also [38])

$$\lim_{N\to\infty}\left[\sum_{k=1}^{N}\frac{H_k-\log k}{k}-\gamma\log N\right]=\frac{1}{2}\left[\varsigma(2)+\gamma^2\right]-\gamma_1$$

and, in the same manner as (i) above, we then obtain the series

$$\sum_{k=1}^{\infty}\left[\frac{H_k-\log k}{k}-\gamma\log\left(1+\frac{1}{k}\right)\right]=\frac{1}{2}\left[\varsigma(2)+\gamma^2\right]-\gamma_1$$

We showed in [38] that

$$\sum_{k=1}^{\infty}\frac{H_k-\gamma-\log k}{k}=\frac{1}{2}\left[\varsigma(2)-\gamma^2\right]-\gamma_1$$

and note that this formula was proposed as a problem by Furdui [39] in 2007. Subtraction of these two formulae simply gives us the well-known result

$$\sum_{k=1}^{\infty}\left[\frac{1}{k}-\log\left(1+\frac{1}{k}\right)\right]=\gamma$$

**(iii) A series for $\varsigma'(0,x)$**

Using the Euler-Maclaurin summation formula, Berndt [3] has shown that (see also [31, p.88])

(4.3) $$\varsigma'(0,x)=\lim_{N\to\infty}\left[-x-N-\log x-\sum_{k=1}^{N}\log(k+x)+(N+x+\tfrac{1}{2})\log(N+x)\right]$$

Using (2.13) we see that

$$\varsigma'(0,1+x)=\varsigma'(0,x)+\log x$$

and hence we obtain

$$\varsigma'(0,1+x)=\lim_{N\to\infty}\left[-x-N-\sum_{k=1}^{N}\log(k+x)+(N+x+\tfrac{1}{2})\log(N+x)\right]$$

which implies that

(4.4) $$\varsigma'(0)=\lim_{N\to\infty}\left[-N-\sum_{k=1}^{N}\log k+(N+\tfrac{1}{2})\log N\right]$$

Letting $x=1$ in (4.3) gives us

$$(4.5) \qquad \varsigma'(0) = \lim_{N \to \infty} \left[ -1 - N - \sum_{k=1}^{N} \log(k+1) + (N + \tfrac{3}{2}) \log(N+1) \right]$$

With $M = N+1$ we see that (4.4) and (4.5) are equivalent.

We then obtain the difference

$$\varsigma'(0,x) - \varsigma'(0) = \lim_{N \to \infty} \left[ -x - \log x - \sum_{k=1}^{N} [\log(k+x) - \log k] + (N + \tfrac{1}{2})[\log(N+x) - \log N] + x\log(N+x) \right]$$

and using (4.1) with Lerch's formula we have

$$= \lim_{N \to \infty} \left[ x \log N - \log x + \sum_{k=1}^{N} \log k - \sum_{k=1}^{N} \log(k+x) \right]$$

Employing Lemma 3.1

$$\log N = -\int_{N}^{N+1} \frac{1}{t} dt + \sum_{k=1}^{N} [\log(k+1) - \log k]$$

we obtain

$$\varsigma'(0,x) - \varsigma'(0) = \lim_{N \to \infty} \left[ -x \int_{N}^{N+1} \frac{1}{t} dt + x \sum_{k=1}^{N} [\log(k+1) - \log k] - \log x + \sum_{k=1}^{N} \log k - \sum_{k=1}^{N} \log(k+x) \right]$$

which may be written as

$$(4.6) \qquad \varsigma'(0,x) - \varsigma'(0) = -\log x + \sum_{k=1}^{N} \Big[ x[\log(k+1) - \log k] + \log k - \log(k+x) \Big]$$

### (iii) The $\eta_k$ constants in terms of the von Mangoldt function

The eta constants $\eta_k$ are defined by reference to the logarithmic derivative of the Riemann zeta function

$$\frac{d}{ds}[\log \varsigma(s)] = \frac{\varsigma'(s)}{\varsigma(s)} = -\frac{1}{s-1} - \sum_{k=1}^{\infty} \eta_{k-1}(s-1)^{k-1} \qquad |s-1| < 3$$

and we also note that this is equivalent to

$$\frac{d}{ds} \log[(s-1)\varsigma(s)] = \frac{\varsigma'(s)}{\varsigma(s)} + \frac{1}{s-1} = -\sum_{k=1}^{\infty} \eta_{k-1}(s-1)^{k-1}$$

and, noting (1.3), $\lim_{s \to 1}[(s-1)\varsigma(s)] = 1$, we obtain upon integration

$$\log[(s-1)\varsigma(s)] = -\sum_{k=1}^{\infty} \frac{\eta_{k-1}}{k}(s-1)^k$$

The $\eta_n$ constants may be written as [16]

$$\eta_n = \frac{(-1)^n}{n!} \lim_{N \to \infty}\left[\sum_{k=1}^{N} \Lambda(k)\frac{\log^n k}{k} - \frac{\log^{n+1} N}{n+1}\right]$$

where $\Lambda$ is the von Mangoldt function [28] such that $\Lambda(k) = \log p$ when $k$ is a power of a prime number $p$ and $\Lambda(k) = 0$ otherwise.

We may therefore employ (3.6) in Lemma 3.1

$$\log^{n+1} N = -\int_{N}^{N+1} \frac{\log^n t}{t} dt + \sum_{k=1}^{N}[\log^{n+1}(k+1) - \log^{n+1} k]$$

to obtain the series representation

(4.7) $$\eta_n = \frac{(-1)^n}{n!} \sum_{k=1}^{\infty}\left[\Lambda(k)\frac{\log^n k}{k} - \frac{1}{n+1}[\log^{n+1}(k+1) - \log^{n+1} k]\right]$$

We recall (3.10)

$$\gamma_n = \sum_{k=1}^{\infty}\left[\frac{\log^n k}{k} - \frac{1}{n+1}[\log^{n+1}(k+1) - \log^{n+1} k]\right]$$

and obtain

(4.8) $$(-1)^n n!\eta_n - \gamma_n = \sum_{k=1}^{\infty}[\Lambda(k) - 1]\frac{\log^n k}{k}$$

We note from the definition of the von Mangoldt function that $\Lambda(k) - 1$ assumes both positive and negative values (thereby enabling the conditional convergence of the above series).

Based on some complex math by Matsuoka [33], Coffey [14] proved that

$$\eta_n = (-1)^{n+1} c_n, \text{ where } c_n > 0$$

and we obtain

$$0 > \sum_{k=1}^{\infty}\left[\Lambda(k)\frac{\log^n k}{k} - \frac{1}{n+1}[\log^{n+1}(k+1) - \log^{n+1} k]\right]$$

### (iv) The $\delta_n$ constants

The $\delta_n$ constants were considered by Sitaramachandrarao [35] in 1986 and by Lehmer [32] in 1988 where they are defined as

$$\varsigma(s) - \frac{1}{s-1} = \sum_{n=0}^{\infty} \frac{(-1)^n \delta_n}{n!} s^n$$

and, since $\varsigma(0) = -\frac{1}{2}$, it is known that $\delta_0 = \frac{1}{2}$ and we have

$$(4.9) \qquad \delta_n = \lim_{N \to \infty} \left[ \sum_{k=1}^{N} \log^n k - \int_{1}^{N} \log^n x \, dx - \frac{1}{2} \log^n N \right]$$

$$= (-1)^n \left[ \varsigma^{(n)}(0) + n! \right]$$

We see that

$$\gamma = \sum_{n=0}^{\infty} \frac{(-1)^n \delta_n}{n!}$$

As pointed out by Apostol [2], the following power series expansion

$$\varsigma(1-s) + \frac{1}{s} = \sum_{n=0}^{\infty} \left[ \frac{\varsigma^{(n)}(0)}{n!} + 1 \right] (1-s)^n$$

converges for $s = 0$ and we therefore have

$$\lim_{n \to \infty} \left[ \frac{\varsigma^{(n)}(0)}{n!} + 1 \right] = 0$$

Therefore we know that for all $\varepsilon > 0$ there exists an $N$ such that $\frac{\varsigma^{(n)}(0)}{n!} + 1 < \varepsilon$ for all $n > N$.

Apostol [2] has calculated the first 18 values of $\frac{\varsigma^{(n)}(0)}{n!}$ and we note that they are all negative and exhibit some small oscillations around the value $-1$.

We also see that

$$\lim_{n \to \infty} \left[ \frac{\varsigma^{(n-1)}(0)}{(n-1)!} + 1 \right] - \lim_{n \to \infty} \left[ \frac{\varsigma^{(n)}(0)}{n!} + 1 \right] = 0$$

and thus we have

$$\lim_{n\to\infty}\frac{1}{n!}\left[n\varsigma^{(n-1)}(0)-\varsigma^{(n)}(0)\right]=0$$

We note that

$$\lim_{s\to 1}\left[\varsigma(s)-\frac{1}{s-1}\right]=\sum_{n=0}^{\infty}\left[\frac{\varsigma^{(n)}(0)}{n!}+1\right]=\gamma$$

A formal derivation of (4.9) follows. Using the Euler-Maclaurin summation formula, Hardy [25, p.333] showed that the Riemann zeta function could be expressed as follows

(4.10) $\qquad \varsigma(s)=\lim_{N\to\infty}\left[\sum_{k=1}^{N}\frac{1}{k^s}-\frac{N^{1-s}}{1-s}-\frac{1}{2}N^{-s}\right]\qquad \operatorname{Re}(s)>-1$

It may immediately be seen that this identity is trivially satisfied for $\operatorname{Re}(s)>1$ because

$$\varsigma(s)=\lim_{N\to\infty}\sum_{k=1}^{N}\frac{1}{k^s}+\lim_{N\to\infty}\left[-\frac{N^{1-s}}{1-s}-\frac{1}{2}N^{-s}\right]$$

and the latter limit is clearly equal to zero.

Differentiating (4.10) results in for $\operatorname{Re}(s)>-1$ (boldly assuming that the operation is valid)

$$\varsigma'(s)=\lim_{N\to\infty}\left[-\sum_{k=1}^{N}\frac{\log k}{k^s}+\frac{N^{1-s}(1-s)\log N-N^{1-s}}{(1-s)^2}+\frac{1}{2}N^{-s}\log N\right]$$

and with $s=0$ we obtain

$$\varsigma'(0)=\lim_{N\to\infty}\left[-\sum_{k=1}^{N}\log k+\left(N+\frac{1}{2}\right)\log N-N\right]$$

Hence, using the Stirling approximation for $n!$ we see that $\varsigma'(0)=-\frac{1}{2}\log(2\pi)$.

With regard to (4.10) we could determine $\varsigma''(0)$

$$\varsigma''(s)=\lim_{N\to\infty}\left[\sum_{k=1}^{N}\frac{\log^2 k}{k^s}+\frac{(1-s)^2[-N^{1-s}(1-s)\log^2 N]+2[N^{1-s}(1-s)\log N-N^{1-s}]}{(1-s)^4}-\frac{1}{2}N^{-s}\log^2 N\right]$$

so that

$$\varsigma''(0)=\lim_{N\to\infty}\left[\sum_{k=1}^{N}\log^2 k-\left(N+\frac{1}{2}\right)\log^2 N+2N\log N-2N\right]$$

In order to simplify the calculations, we write

$$\varsigma(s) = \lim_{N\to\infty}\left[\sum_{k=1}^{N}\frac{1}{k^s} - \frac{N^{1-s}}{1-s} - \frac{1}{2}N^{-s}\right] = \lim_{N\to\infty}\left[\sum_{k=1}^{N}\frac{1}{k^s} - \frac{N^{1-s}-1}{1-s} - \frac{1}{1-s} - \frac{1}{2}N^{-s}\right]$$

and formal differentiation gives us (boldly assuming that the operation is valid)

$$\varsigma^{(n)}(s) = \lim_{N\to\infty}\left[(-1)^n\sum_{k=1}^{N}\frac{\log^n k}{k^s} - f^{(n)}(s) - \frac{n!}{(1-s)^{n+1}} - \frac{1}{2}(-1)^n N^{-s}\log^n N\right]$$

where we have denoted $f(s)$ as

$$f(s) = \frac{N^{1-s}-1}{s-1}$$

We can represent $f(s)$ by the following integral

$$f(s) = \frac{N^{1-s}-1}{s-1} = -\int_{1}^{N} x^{-s}\,dx$$

so that

(4.11) $\qquad f^{(n)}(s) = -(-1)^n \int_{1}^{N} x^{-s}\log^n x\,dx$

and thus

$$f^{(n)}(0) = -(-1)^n \int_{1}^{N}\log^n x\,dx$$

Therefore, with $s=0$ we obtain [13]

(4.12) $\qquad (-1)^n\left[\varsigma^{(n)}(0) + n!\right] = \lim_{N\to\infty}\left[\sum_{k=1}^{N}\log^n k - \int_{1}^{N}\log^n x\,dx - \frac{1}{2}\log^n N\right]$

With the substitution $y = \log(u+k)$ we get

$$\int \log^{p+1}(u+k)\,du = \int y^{p+1}e^y\,dy = (-1)^{p+1}(p+1)!\, e^y \sum_{m=0}^{p+1}(-1)^m \frac{y^m}{m!}$$

The latter integral may be easily obtained by parametric differentiation.

Therefore we have

$$\int \log^{p+1}(u+k)\,du = (-1)^{p+1}(p+1)!\sum_{m=0}^{p+1}\frac{(-1)^m}{m!}(u+k)\log^m(u+k)$$

and in particular we have

$$\int_0^u \log^n x\, dx = (-1)^n n!\, u \sum_{j=0}^n (-1)^j \frac{\log^j u}{j!}$$

and

$$\int_0^1 \log^n x\, dx = (-1)^n n!$$

and thus

$$\int_1^u \log^n x\, dx = -(-1)^n n! + (-1)^n n!\, u \sum_{j=0}^n (-1)^j \frac{\log^j u}{j!}$$

In passing we note that

$$\int_1^N \log^n x\, dx = (-1)^n \int_1^N \log^n (1/x)\, dx$$

$$= -(-1)^n \int_1^{1/N} \frac{\log^n t}{t^2}\, dt$$

and Choudhury [13] has noted the evaluation of this integral.

We see that

$$(-1)^n \varsigma^{(n)}(0) = \lim_{N \to \infty} \left[ \sum_{k=1}^N \log^n k - (-1)^n n!\, N \sum_{j=0}^n (-1)^j \frac{\log^j N}{j!} - \frac{1}{2} \log^n N \right]$$

$\square$

We see that

$$\varsigma(s) - \frac{1}{s-1} = \lim_{N \to \infty} \left[ \sum_{k=1}^N \frac{1}{k^s} - \frac{N^{1-s} - 1}{1-s} - \frac{1}{2} N^{-s} \right]$$

and formal differentiation gives us (boldly assuming that the operation is valid)

$$D^n \left[ \varsigma(s) - \frac{1}{s-1} \right] = \lim_{N \to \infty} \left[ (-1)^n \sum_{k=1}^N \frac{\log^n k}{k^s} - f^{(n)}(s) - \frac{1}{2} (-1)^n N^{-s} \log^n N \right]$$

$$D^n \left[ \varsigma(s) - \frac{1}{s-1} \right]_{s=1} = \lim_{N \to \infty} \left[ (-1)^n \sum_{k=1}^N \frac{\log^n k}{k} - f^{(n)}(1) - \frac{1}{2} (-1)^n N^{-1} \log^n N \right]$$

We see from (4.11) that

$$f^{(n)}(1) = (-1)^{n+1} \int_1^N \frac{\log^n t}{t}\, dt = (-1)^{n+1} \frac{\log^{n+1} N}{n+1}$$

which results in

$$D^n\left[\varsigma(s)-\frac{1}{s-1}\right]_{s=1} = \lim_{N\to\infty}\left[(-1)^n\sum_{k=1}^{N}\frac{\log^n k}{k} - (-1)^n\frac{\log^{n+1} N}{n+1} - \frac{1}{2}(-1)^n N^{-1}\log^n N\right]$$

Hence we obtain

$$\gamma_n = \lim_{N\to\infty}\left[\sum_{k=1}^{N}\frac{\log^n k}{k} - \frac{\log^{n+1} N}{n+1} - \frac{1}{2}N^{-1}\log^n N\right]$$

$$= \lim_{N\to\infty}\left[\sum_{k=1}^{N}\frac{\log^n k}{k} - \frac{\log^{n+1} N}{n+1}\right]$$

which accords with (3.1).

Similar exercises could be undertaken in respect of the known limit formulae for the Glaisher-Kinkelin constants [12] (albeit involving more complexity).

## 5. Some relations involving the derivatives of the Hurwitz zeta function $\varsigma^{(k)}(0,x)$

**Proposition 5.1**

$$\varsigma''(0,x) - \varsigma''(0) = \log^2 x + \sum_{n=1}^{\infty}\left[\log^2(n+x) - \log^2 n - x[\log^2(n+1) - \log^2 n]\right]$$

**Proof**

We define $g_k(x)$ as the function such that $g_k''(x) := \varsigma^{(k)}(2,x)$ where $k \geq 1$ is an integer and, for convenience, we specify the initial conditions $g_k'(1) = g_k(1) = 0$.

We have in particular using (2.1)

$$g_1''(x) := \varsigma'(2,x) = -\sum_{n=0}^{\infty}\frac{\log(n+x)}{(n+x)^2}$$

Since

$$\int\frac{\log(n+x)}{(n+x)^2}dx = -\left[\frac{\log(n+x)}{n+x} + \frac{1}{n+x}\right]$$

integration results in

(5.1) $$g_1'(x) = \sum_{n=0}^{\infty}\left[\frac{\log(n+x)}{n+x} - \frac{\log(n+1)}{n+1} + \frac{1}{n+x} - \frac{1}{n+1}\right]$$

since we have specified ab initio that $g_1'(1) = 0$.

Referring to (2.7) we see that

$$\gamma_1(x) - \gamma_1 = \sum_{n=0}^{\infty}\left[\frac{\log(n+x)}{n+x} - \frac{\log(n+1)}{n+1}\right]$$

and hence we have

(5.2) $$g_1'(x) = \gamma_1(x) - \gamma_1 - [\psi(x) + \gamma]$$

Another integration gives us

(5.3) $$g_1(x) = \int_1^x \gamma_1(x)dx - (x-1)\gamma_1 - [\log \Gamma(x) + (x-1)\gamma]$$

Referring to (2.16) we see that

$$\int_1^x \gamma_1(t)\, dt = \frac{1}{2}\left[\varsigma''(0,x) - \varsigma''(0)\right]$$

which gives us

(5.4) $$g_1(x) = \frac{1}{2}\left[\varsigma''(0,x) - \varsigma''(0)\right] - (x-1)\gamma_1 - [\log \Gamma(x) + (x-1)\gamma]$$

Alternatively, integrating (5.1) gives us

(5.5) $$g_1(x) = \frac{1}{2}\sum_{n=0}^{\infty}\left[\log^2(n+x) - \log^2(n+1) - 2(x-1)\frac{\log(n+1)}{n+1}\right] - [\log \Gamma(x) + (x-1)\gamma]$$

Equating (5.4) and (5.5) we obtain

(5.6) $$\varsigma''(0,x) - \varsigma''(0) - 2(x-1)\gamma_1 = \sum_{n=0}^{\infty}\left[\log^2(n+x) - \log^2(n+1) - 2(x-1)\frac{\log(n+1)}{n+1}\right]$$

Letting $x \to 1+x$ in (5.6) results in

$$\varsigma''(0,1+x) - \varsigma''(0) - 2\gamma_1 x = \sum_{n=0}^{\infty}\left[\log^2(n+1+x) - \log^2(n+1) - 2x\frac{\log(n+1)}{n+1}\right]$$

and, since $\varsigma''(0,1+x) = \varsigma''(0,x) - \log^2 x$, we obtain

(5.7) $$\varsigma''(0,x) - \varsigma''(0) - 2\gamma_1 x = \log^2 x + \sum_{n=1}^{\infty}\left[\log^2(n+x) - \log^2 n - 2x\frac{\log n}{n}\right]$$

This is in agreement with Deninger's result [22] where he defined $R(x) = -\varsigma''(0,x)$ (it should be noted that he employs a slightly different definition of $\gamma_1$).

With $x = 1$ we have

(5.8) $$\gamma_1 = \sum_{n=1}^{\infty}\left[\frac{\log n}{n} - \frac{1}{2}[\log^2(n+1) - \log^2 n]\right]$$

and this concurs with (3.13).

With $x = 2$ we have

$$\gamma_1 = -\frac{1}{4}\log^2 2 + \sum_{n=1}^{\infty}\left[\frac{\log n}{n} - \frac{1}{4}[\log^2(n+2) - \log^2 n]\right]$$

Using (5.7) and (5.8) we obtain

(5.9) $$\varsigma''(0, x) - \varsigma''(0) = \log^2 x + \sum_{n=1}^{\infty}\left[\log^2(n+x) - \log^2 n - x[\log^2(n+1) - \log^2 n]\right]$$

We note Kanemitsu and Tsukada [31, p.90] report that

(5.10) $$\varsigma''(0, x) - \varsigma''(0) = \log^2 x + \lim_{N \to \infty}\left(\sum_{n=1}^{N}\left[\log^2(n+x) - \log^2 n\right] - x\log^2 N\right)$$

Lemma 3.1 gives us

$$\log^2 N = -\int_N^{N+1} \frac{\log^n t}{t} dt + \sum_{n=1}^{N}[\log^2(n+1) - \log^2 n]$$

and hence we obtain

$$\lim_{N \to \infty}\left(\sum_{n=1}^{N}\left[\log^2(n+x) - \log^2 n\right] - x\log^2 N\right)$$

$$= \lim_{N \to \infty}\left(\sum_{n=1}^{N}\left[\log^2(n+x) - \log^2 n - x[\log^2(n+1) - \log^2 n]\right]\right)$$

Then taking the limit $N \to \infty$ results in (5.9).

The same source [31, p.91] reports that

$$\varsigma''(0, x) - \varsigma''(0) \stackrel{??}{=} 2\gamma_1 \log x - \log^2 x + \sum_{n=1}^{N}\left[\log^2(n+x) - \log^2 n - 2\frac{\log n}{n}\right] + O\left(\frac{\log^2 N}{N}\right)$$

but there are two typos in this formula and, as shown below, the term $\log x$ should be replaced by $x$ and the term $\frac{\log n}{n}$ should be multiplied by $x$ (this formula was also incorrectly quoted in [30]). The corrected formula is displayed below.

$$(5.11) \quad \varsigma''(0,x) - \varsigma''(0) = 2\gamma_1 x - \log^2 x + \sum_{k=1}^{N} \left[ \log^2(k+x) - \log^2 k - 2x \frac{\log k}{k} \right] + O\left( \frac{\log^2 N}{N} \right)$$

Taking the limit $N \to \infty$ results in (5.7).

□

Integration of (5.9) over $[0,u]$ gives us

$$\int_0^u \varsigma''(0,x)\,dx - \varsigma''(0)u = u[\log^2 u - 2\log u + 2]$$

$$+ \sum_{n=1}^{\infty} \left[ (n+u)\left\{\log^2(n+u) - 2\log(n+u) + 2\right\} - n\left\{\log^2 n - 2\log n + 2\right\} - u\log^2 n - \frac{1}{2}u^2[\log^2(n+1) - \log^2 n] \right]$$

As noted by Deninger [22] we have

$$(5.12) \quad \int_0^1 \varsigma''(0,x)\,dx = 0$$

and we therefore have

(5.13)
$$\varsigma''(0) = -2 - \sum_{n=1}^{\infty} \left[ (n+1)\left\{\log^2(n+1) - 2\log(n+1) + 2\right\} - n\left\{\log^2 n - 2\log n + 2\right\} - \log^2 n - \frac{1}{2}[\log^2(n+1) - \log^2 n] \right]$$

Choudhury [13] reports that

$$(5.14) \quad \varsigma''(0) = -2 + o(1)$$

□

We may approach Proposition 5.1 in a slightly different direction as set out below.

We note from (1.1) that

$$\frac{\partial}{\partial x}\varsigma(s,x) = -s\varsigma(s+1,x)$$

and a few partial derivatives are set out below.

$$\frac{\partial^2}{\partial x^2}\varsigma(s,x) = s(s+1)\varsigma(s+2,x)$$

$$\frac{\partial}{\partial s}\frac{\partial^2}{\partial x^2}\varsigma(s,x) = s(s+1)\varsigma'(s+2,x)+(2s+1)\varsigma(s+2,x)$$

$$\frac{\partial^2}{\partial s^2}\frac{\partial^2}{\partial x^2}\varsigma(s,x) = s(s+1)\varsigma''(s+2,x)+2(2s+1)\varsigma'(s+2,x)+2\varsigma(s+2,x)$$

With $s=0$ we have

$$\frac{\partial^2}{\partial x^2}\varsigma''(0,x) = 2\varsigma'(2,x)+2\varsigma(2,x)$$

which may be expressed as

$$\frac{\partial^2}{\partial x^2}\varsigma''(0,x) = 2\varsigma'(2,x)+2\frac{d^2}{dx^2}\log\Gamma(x)$$

Using the previous notation $g_1''(x) := \varsigma'(2,x) = -\sum_{n=0}^{\infty}\frac{\log(n+x)}{(n+x)^2}$ we have

$$\frac{d^2}{dx^2}[\varsigma''(0,x) - 2g_1(x) - 2\log\Gamma(x)] = 0$$

and integration results in

$$\varsigma''(0,x) - 2g_1(x) - 2\log\Gamma(x) = ax+b$$

The integration constants are easily determined as follows. We have the derivative

$$\frac{d}{dx}\varsigma''(0,x) - 2g_1'(x) - 2\psi(x) = a$$

and we see from (2.15) that

$$\frac{d}{dx}\varsigma''(0,x) = 2\gamma_1(x)$$

and thus

$$2\gamma_1(x) - 2g_1'(x) - 2\psi(x) = a$$

With the initial conditions $g_1'(1) = g_1(1) = 0$, a little algebra shows that

$$g_1(x) = \frac{1}{2}[\varsigma''(0,x) - \varsigma''(0)] - (x-1)\gamma_1 - [\log\Gamma(x) - (x-1)\gamma]$$

**Proposition 5.2**

(5.14) $\quad \varsigma^{(3)}(0,x) - \varsigma^{(3)}(0) = -\log^3 x - \sum_{n=1}^{\infty} \left[ \log^3(n+x) - \log^3 n - x[\log^3(n+1) - \log^3 n] \right]$

**Proof**

With reference to the previous proposition, we now consider the case

$$g_2''(x) := \varsigma''(2,x) = \sum_{n=0}^{\infty} \frac{\log^2(n+x)}{(n+x)^2}$$

Since

$$\int \frac{\log^2(n+x)}{(n+x)^2} dx = -2 \left[ \frac{1}{2} \frac{\log^2(n+x)}{n+x} + \frac{\log(n+x)}{n+x} + \frac{1}{n+x} \right]$$

integration results in

$$g_2'(x) = -2 \sum_{n=0}^{\infty} \left[ \frac{1}{2} \frac{\log^2(n+x)}{n+x} - \frac{1}{2} \frac{\log^2(n+1)}{n+1} + \frac{\log(n+x)}{n+x} - \frac{\log(n+1)}{n+1} + \frac{1}{n+x} - \frac{1}{n+1} \right]$$

which, having regard to (5.1), we may express as

$$g_2'(x) = -\sum_{n=0}^{\infty} \left[ \frac{\log^2(n+x)}{n+x} - \frac{\log^2(n+1)}{n+1} \right] - 2g_1'(x)$$

Referring to (2.7) we see that

$$\gamma_2(x) - \gamma_2 = \sum_{n=0}^{\infty} \left[ \frac{\log^2(n+x)}{n+x} - \frac{\log^2(n+1)}{n+1} \right]$$

and thus we have

$$g_2'(x) = -[\gamma_2(x) - \gamma_2] - 2g_1'(x)$$

Another integration gives us

$$g_2(x) = \int_1^x \gamma_2(x) dx - (x-1)\gamma_2 - 2g_1(x)$$

since $g_2(1) = 0$. We note from (2.16) that

$$\int_1^x \gamma_2(t) dt = -\frac{1}{3} \left[ \varsigma^{(3)}(0,x) - \varsigma^{(3)}(0) \right]$$

and we therefore obtain

(5.15) $$g_2(x) = -\frac{1}{3}\left[\varsigma^{(3)}(0,x) - \varsigma^{(3)}(0)\right] - (x-1)\gamma_2 - 2g_1(x)$$

Alternatively, using the form

$$g_2'(x) = -\sum_{n=0}^{\infty}\left[\frac{\log^2(n+x)}{n+x} - \frac{\log^2(n+1)}{n+1}\right] - 2g_1'(x)$$

we integrate to obtain

(5.16) $$g_2(x) = \frac{1}{3}\sum_{n=0}^{\infty}\left[\log^3(n+x) - \log^3(n+1) - 3(x-1)\frac{\log^2(n+1)}{n+1}\right] - 2g_1(x)$$

Equating (5.15) and (5.16) we obtain

(5.17) $$\varsigma^{(3)}(0,x) - \varsigma^{(3)}(0) + 3(x-1)\gamma_2 = -\sum_{n=0}^{\infty}\left[\log^3(n+x) - \log^3(n+1) - 3(x-1)\frac{\log^2(n+1)}{n+1}\right]$$

Letting $x \to 1+x$ in (5.17) and using

$$\varsigma^{(n)}(0, 1+x) = \varsigma^{(n)}(0, x) + (-1)^{n+1}\log^n x$$

results in

$$\varsigma^{(3)}(0,x) - \varsigma^{(3)}(0) + 3x\gamma_2 = -\log^3 x - \sum_{n=0}^{\infty}\left[\log^3(n+1+x) - \log^3(n+1) - 3x\frac{\log^2(n+1)}{n+1}\right]$$

With $x=1$ we obtain

$$\gamma_2 = \sum_{n=0}^{\infty}\left\{\frac{\log^2(n+1)}{n+1} - \frac{1}{3}[\log^3(n+2) - \log^3(n+1)]\right\}$$

in agreement with (3.14). We then see that

$$\varsigma^{(3)}(0,x) - \varsigma^{(3)}(0) = -\log^3 x - \sum_{n=0}^{\infty}\left[\log^3(n+1+x) - \log^3(n+1) - x[\log^3(n+2) - \log^3(n+1)]\right]$$

which may be expressed by re-indexing as

$$\varsigma^{(3)}(0,x) - \varsigma^{(3)}(0) = -\log^3 x - \sum_{n=1}^{\infty}\left[\log^3(n+x) - \log^3 n - x[\log^3(n+1) - \log^3 n]\right]$$

and this may be compared with

$$\varsigma''(0,x)-\varsigma''(0)=\log^2 x+\sum_{n=1}^{\infty}\left[\log^2(n+x)-\log^2 n-x[\log^2(n+1)-\log^2 n]\right]$$

As noted above in (4.2) we have

$$\log\Gamma(x+1)=\sum_{n=1}^{\infty}\left[x\log\left(1+\frac{1}{n}\right)-\log\left(1+\frac{x}{n}\right)\right]$$

which may be expressed in a similar format

$$\log\Gamma(x)=-\log x-\sum_{n=1}^{\infty}\left[\log(n+x)-\log n-x[\log(n+1)-\log n]\right]$$

Using Lerch's formula $\log\Gamma(x)=\varsigma'(0,x)-\varsigma'(0)$ this becomes

$$\varsigma'(0,x)-\varsigma'(0)=-\log x-\sum_{n=1}^{\infty}\left[\log(n+x)-\log n-x[\log(n+1)-\log n]\right]$$

These three results lead one to suspect that we have the following generalisation.

**Proposition 5.3**

(5.18) $\quad (-1)^{k+1}[\varsigma^{(k+1)}(0,x)-\varsigma^{(k+1)}(0)]$

$$=\log^{k+1}x+\sum_{n=0}^{\infty}\left[\log^{k+1}(n+1+x)-\log^{k+1}(n+1)-x[\log^{k+1}(n+2)-\log^{k+1}(n+1)]\right]$$

where we have started the summation at $n=0$.

This may also be expressed by re-indexing as

(5.19)
$$(-1)^{k+1}[\varsigma^{(k+1)}(0,x)-\varsigma^{(k+1)}(0)]=\log^{k+1}x+\sum_{n=1}^{\infty}\left[\log^{k+1}(n+x)-\log^{k+1}n-x[\log^{k+1}(n+1)-\log^{k+1}n]\right]$$

**Proof**

Combining (5.18) with (2.16)

$$\int_{1}^{u}\gamma_k(x)\,dx=\frac{(-1)^{k+1}}{k+1}\left[\varsigma^{(k+1)}(0,u)-\varsigma^{(k+1)}(0)\right]$$

gives us

$$(k+1)\int_{1}^{u}\gamma_k(x)\,dx=\log^{k+1}x+\sum_{n=0}^{\infty}\left[\log^{k+1}(n+1+x)-\log^{k+1}(n+1)-x[\log^{k+1}(n+2)-\log^{k+1}(n+1)]\right]$$

Differentiation results in

$$(5.20) \quad \gamma_k(x) = \frac{\log^k x}{x} + \sum_{n=0}^{\infty}\left[\frac{\log^k(n+1+x)}{n+1+x} - \frac{1}{k+1}[\log^{k+1}(n+2) - \log^{k+1}(n+1)]\right]$$

which, prima facie, appears quite different from (3.9) which is reproduced below

$$\gamma_k(x) = -\frac{\log^{k+1} x}{k+1} + \sum_{n=0}^{\infty}\left[\frac{\log^k(n+x)}{n+x} - \frac{1}{k+1}[\log^{k+1}(n+1+x) - \log^{k+1}(n+x)]\right]$$

Subtracting these two equations gives us

$$\Delta := \frac{\log^k x}{x} + \sum_{n=0}^{\infty}\left[\frac{\log^k(n+1+x)}{n+1+x} - \frac{1}{k+1}[\log^{k+1}(n+2) - \log^{k+1}(n+1)]\right]$$

$$+ \frac{\log^{k+1} x}{k+1} - \sum_{n=0}^{\infty}\left[\frac{\log^k(n+x)}{n+x} - \frac{1}{k+1}[\log^{k+1}(n+1+x) - \log^{k+1}(n+x)]\right]$$

Using the finite telescoping sum $\sum_{k=0}^{N}[a_{k+1} - a_k] = a_{N+1} - a_0$ we obtain

$$\sum_{n=0}^{N}[\log^{k+1}(n+1+x) - \log^{k+1}(n+x)] = \log^{k+1}(N+1+x) - \log^{k+1} x$$

$$\sum_{n=0}^{N}[\log^{k+1}(n+1+1) - \log^{k+1}(n+1)] = \log^{k+1}(N+2)$$

and we then have

$$\Delta = \frac{\log^k x}{x} + \frac{1}{k+1}[\log^{k+1}(N+1+x) - \log^{k+1} x] + \sum_{n=0}^{N}\left[\frac{\log^k(n+1+x)}{n+1+x} - \frac{\log^k(n+x)}{n+x}\right]$$

$$+ \frac{\log^{k+1} x}{k+1} - \frac{1}{k+1}\log^{k+1}(N+2)$$

Therefore, as $N \to \infty$ we obtain

$$\Delta = \frac{\log^k x}{x} + \sum_{n=0}^{\infty}\left[\frac{\log^k(n+1+x)}{n+1+x} - \frac{\log^k(n+x)}{n+x}\right]$$

$$= \frac{\log^k x}{x} + \gamma_k(1+x) - \gamma_k(x)$$

and reference to (2.10) shows that $\Delta = 0$.

$\square$

Since $\gamma_n(1+x) - \gamma_n(x) = -\dfrac{\log^n x}{x}$ we see from (5.20) that

$$\gamma_k(1+x) = \sum_{n=0}^{\infty}\left[\frac{\log^k(n+1+x)}{n+1+x} - \frac{1}{k+1}[\log^{k+1}(n+2) - \log^{k+1}(n+1)]\right]$$

and this implies that

(5.21) $$\gamma_k(x) = \sum_{n=0}^{\infty}\left[\frac{\log^k(n+x)}{n+x} - \frac{1}{k+1}[\log^{k+1}(n+2) - \log^{k+1}(n+1)]\right]$$

With $x=1$ we have

$$\gamma_k = \sum_{n=0}^{\infty}\left[\frac{\log^k(n+1)}{n+1} - \frac{1}{k+1}[\log^{k+1}(n+2) - \log^{k+1}(n+1)]\right]$$

$\square$

Differentiation of (5.21) results in

$$\gamma_k'(x) = \sum_{n=0}^{\infty}\left[\frac{k\log^{k-1}(n+x) - \log^k(n+x)}{(n+x)^2}\right]$$

which may be written as

$$\gamma_k'(x) = (-1)^{k-1}[k\varsigma^{(k-1)}(2,x) + \varsigma^{(k)}(2,x)]$$

It should be noted that this may be obtained more directly by differentiating (2.7).

$$\gamma_k(x) - \gamma_k = \sum_{n=0}^{\infty}\left[\frac{\log^k(n+x)}{n+x} - \frac{\log^k(n+1)}{n+1}\right]$$

$\square$

With reference to (5.18) we have

$$S := \sum_{n=0}^{\infty}\left[\log^{k+1}(n+1+x) - \log^{k+1}(n+1) - x[\log^{k+1}(n+2) - \log^{k+1}(n+1)]\right]$$

$$= (k+1)\sum_{n=0}^{\infty}\int_1^{1+x}\left[\frac{\log^k(n+u)}{n+u} - \frac{1}{k+1}[\log^{k+1}(n+2) - \log^{k+1}(n+1)]\right]du$$

We next assume that we may change the order of summation and integration

$$= (k+1) \int_{1}^{1+x} \sum_{n=0}^{\infty} \left[ \frac{\log^{k}(n+u)}{n+u} - \frac{1}{k+1} [\log^{k+1}(n+2) - \log^{k+1}(n+1)] \right] du$$

and referring to (5.21) we obtain

$$S = (k+1) \int_{1}^{1+x} \gamma_{k}(u)\, du$$

From (2.16) we see that

$$\int_{1}^{1+x} \gamma_{k}(u)\, du = \frac{(-1)^{k+1}}{k+1} [\varsigma^{(k+1)}(0, 1+x) - \varsigma^{(k+1)}(0)]$$

and thus

$$S = (-1)^{k+1} [\varsigma^{(k+1)}(0, 1+x) - \varsigma^{(k+1)}(0)]$$

Using (2.13)

$$\varsigma^{(k+1)}(0, 1+x) - \varsigma^{(k+1)}(0, x) = (-1)^{k} \log^{k+1} x$$

we obtain

$$S = (-1)^{k+1} [\varsigma^{(k+1)}(0, x) - \varsigma^{(k+1)}(0)] - \log^{k+1} x$$

This then proves (5.18).

□

Integration of (5.21) gives us

$$\int_{1}^{u} \gamma_{k}(x)\, dx = \frac{1}{k+1} \sum_{n=0}^{\infty} \left[ \log^{k+1}(n+u) - \log^{k+1}(n+1) - (u-1)[\log^{k+1}(n+2) - \log^{k+1}(n+1)] \right]$$

and (2.16) results in

(5.22)  $(-1)^{k+1} [\varsigma^{(k+1)}(0, x) - \varsigma^{(k+1)}(0)]$

$$= \log^{k+1} x + \sum_{n=0}^{\infty} \left[ \log^{k+1}(n+1+x) - \log^{k+1}(n+1) - x[\log^{k+1}(n+2) - \log^{k+1}(n+1)] \right]$$

□

We see from (5.18) that

$$\lim_{x \to 0} [\varsigma^{(k+1)}(0, x) - \log^{k+1} x] = (-1)^{k+1} \varsigma^{(k+1)}(0)$$

and this is immediately seen to be the case using (2.13)

$$\varsigma^{(k)}(0,1+x) = \varsigma^{(k)}(0,x) + (-1)^{k+1}\log^k x$$

The case $k=0$ is easily checked using Lerch's formula (2.17).

## 6. Some connections with Dilcher's generalised gamma functions

Dilcher [24] defined generalised gamma functions by

$$\Gamma_k(x) = \lim_{n\to\infty} \frac{\exp\left[\dfrac{x}{k+1}\log^{k+1} n\right]\prod_{j=1}^{n}\exp\left[\dfrac{1}{k+1}\log^{k+1} j\right]}{\prod_{j=0}^{n}\exp\left[\dfrac{1}{k+1}\log^{k+1}(j+x)\right]}$$

where we note that $\Gamma_0(x)$ corresponds with the classical gamma function

$$\Gamma_0(x) = \lim_{n\to\infty}\frac{n^x n!}{x(x+1)\cdots(x+n)} = \Gamma(x)$$

and we have

$$\Gamma_k(1) = 1$$

$$\Gamma_k(x+1) = \exp\left[\frac{1}{k+1}\log^{k+1} x\right]\Gamma_k(x)$$

Dilcher [24] showed that

(6.1) $\qquad \log\Gamma_1(x+1) + \gamma_1 x = \sum_{n=1}^{\infty}\dfrac{(-1)^n}{n+1}[H_n\varsigma(n+1) + \varsigma'(n+1)]x^{n+1}$

By definition we have

$$\log\Gamma_1(x) = \lim_{n\to\infty}\left[\frac{x}{2}\log^2 n + \frac{1}{2}\sum_{j=1}^{n}\log^2 j - \frac{1}{2}\sum_{j=0}^{n}\log^2(j+x)\right]$$

and we have

$$\log\Gamma_1(x+1) = \frac{1}{2}\log^2 x + \log\Gamma_1(x)$$

Therefore we have

$$\log\Gamma_1(x+1) + \gamma_1 x = \lim_{n\to\infty}\left[\frac{x}{2}\log^2 n + \frac{1}{2}\sum_{j=1}^{n}\log^2 j - \frac{1}{2}\sum_{j=1}^{n}\log^2(j+x)\right] + \lim_{n\to\infty}\sum_{j=1}^{n}\left[x\frac{\log j}{j} - \frac{x}{2}\log^2 n\right]$$

$$= \lim_{n \to \infty} \left[ x \sum_{j=1}^{n} \frac{\log j}{j} + \frac{1}{2} \sum_{j=1}^{n} \log^2 j - \frac{1}{2} \sum_{j=1}^{n} \log^2(j+x) \right]$$

We showed in [18] that

(6.2) $$\sum_{n=1}^{\infty} \frac{(-1)^n}{n+1}[H_n \varsigma(n+1) + \varsigma'(n+1)] x^{n+1} = \frac{1}{2} \sum_{n=1}^{\infty} \left[ 2x \frac{\log n}{n} + \log^2 n - \log^2(n+x) \right]$$

and hence, referring to (5.7) above, we deduce that

(6.3) $$\varsigma''(0,x) - \varsigma''(0) - 2\gamma_1 x = \log^2 x - 2 \sum_{n=1}^{\infty} \frac{(-1)^n}{n+1}[H_n \varsigma(n+1) + \varsigma'(n+1)] x^{n+1}$$

With $x=1$ we obtain (see [17] and [18])

(6.4) $$\gamma_1 = \sum_{n=1}^{\infty} \frac{(-1)^n}{n+1}[H_n \varsigma(n+1) + \varsigma'(n+1)]$$

It was also shown by Dilcher [24] that

$$\frac{1}{\Gamma_k(x)} = e^{\gamma_k x} \exp\left[\frac{x}{k+1} \log^{k+1} x\right] \prod_{n=1}^{\infty} \exp\left[-\frac{x}{n} \log^k n\right] \exp\left[\frac{1}{k+1}[\log^{k+1}(n+x) - \log^{k+1} n]\right]$$

whereupon it directly follows that

$$\log \Gamma_k(x+1) + \gamma_k x = \lim_{n \to \infty} \left[ x \sum_{j=1}^{n} \frac{\log^k j}{j} + \frac{1}{k+1} \sum_{j=1}^{n} \log^k j - \frac{1}{k+1} \sum_{j=1}^{n} \log^k(j+x) \right]$$

and thus

$$\log \Gamma_k(x+1) + \gamma_k x = \sum_{j=1}^{\infty} \left[ x \frac{\log^k j}{j} - \frac{1}{k+1}[\log^k(j+x) - \log^k j] \right]$$

Since $\Gamma_k(2) = 1$ we see that

$$\gamma_k = \sum_{j=1}^{\infty} \left[ \frac{\log^k j}{j} - \frac{1}{k+1}[\log^k(j+1) - \log^k j] \right]$$

in accordance with (3.10).

Chakraborty, Kanemitsu and Li [11] have hinted at a relationship between the Stieltjes constants and the Barnes' multiple gamma functions in 2010 but I have not been able to find any follow up research on that aspect.

## 7. Open access to our own work

This paper contains references to various other papers and, rather surprisingly, most of them are currently freely available on the internet. Surely now is the time that <u>all</u> of <u>our</u> work should be freely accessible by <u>all</u>. The mathematics community should lead the way on this by publishing <u>everything</u> on arXiv, or in an equivalent open access repository. We think it, we write it, so why hide it? You know it makes sense.

[24] K. Dilcher, On generalized gamma functions related to the Laurent coefficients of the Riemann zeta function. Aeq. Math. 48 (1994) 55-85.

[25] G.H. Hardy, Divergent Series. Chelsea Publishing Company, New York, 1991.

[26] H. Hasse, Ein Summierungsverfahren für Die Riemannsche $\varsigma$ - Reithe. Math.Z.32, 458-464, 1930.
http://dz-srv1.sub.uni-goettingen.de/sub/digbib/loader?ht=VIEW&did=D23956&p=462

[27] C. Hermite, Correspondance d'Hermite et de Stieltjes. Gauthier-Villars, Paris, 1905.
http://ebooks.library.cornell.edu/m/math/browse/author/a.php

[28] A. Ivić, The Riemann Zeta- Function: Theory and Applications. Dover Publications Inc, 2003.

[29] J. L. W. V. Jensen, Sur la fonction ζ(s) de Riemann, Comptes-rendus hebdomadaires des seances de l'Academie des sciences, tome 104, pp. 1156–1159 (1887).

[30] S. Kanemitsu, Y. Tanigawa, H. Tsukada and M. Yoshimoto. Contributions to the theory of the Hurwitz zeta function.
Hardy-Ramanujan Journal, Hardy-Ramanujan Society, 2007, 30, pp.31-55.
https://hal.archives-ouvertes.fr/hal-01112081/document

[31] S. Kanemitsu and H. Tsukada, Vistas of special functions. World Scientific Publishing Co. Pte. Ltd., 2007.

[32] D.H. Lehmer, The sum of like powers of the zeros of the Riemann zeta function. Math. Comp. 50, 265–27,1988.

[33] Y. Matsuoka, Generalized Euler constants associated with the Riemann zeta function. In Number theory and combinatorics (eds J. Akiyama et al.). Singapore: World Scientific. 1985.

[34] T.J. Osler, An Introduction to the Zeta Function. 2004
https://academics.rowan.edu/csm/departments/math/facultystaff/faculty/osler/Introduction-to-Complex-Analysis/ComplexAnalysisDownloads.html

[35] R. Sitaramachandrarao, Maclaurin Coefficients of the Riemann Zeta Function. Abstracts Amer. Math. Soc. **7**, 280, 1986.

[36] H.M. Srivastava and J. Choi, Series Associated with the Zeta and Related Functions. Kluwer Academic Publishers, Dordrecht, the Netherlands, 2001.

[37] E.C. Titchmarsh, The Zeta-Function of Riemann. Oxford University (Clarendon) Press, Oxford, London and New York, 1951; Second Ed. (Revised by D.R. Heath- Brown), 1986.

[38] D.F. Connon, On an integral involving the digamma function. 2012.
arXiv:1212.1432 [pdf]

[39] O. Furdui, College Math. Journal, 38, No.1, 61, 2007
45

dconnon@btopenworld.com

Wessex House,
Devizes Road,
Upavon,
Pewsey,
Wiltshire SN9 6DL